\documentclass[a4paper,11pt,english]{smfart}

\usepackage{amssymb}
\usepackage{url}
\usepackage[T1]{fontenc}
\RequirePackage{calrsfs}
\DeclareSymbolFont{rsfscript}{OMS}{rsfs}{m}{b}
\DeclareSymbolFontAlphabet{\mathrsfs}{rsfscript}
\usepackage[utopia,expert]{mathdesign}
\usepackage{lscape}
\usepackage{array, boldline, makecell, booktabs}

\usepackage{enumitem}

\usepackage{todonotes}

\usepackage[leftbars]{changebar}

\usepackage{palatino}
\usepackage{rotating}
\usepackage{graphicx}

\usepackage{tikz-cd}

\usepackage{enumerate}

\usepackage{amscd}

\usepackage{color}
\definecolor{shadecolor}{gray}{0.90}

\input xypic
\xyoption{all}
\xyoption{arc}

\makeindex

\def\equat{\refstepcounter{theo}\begin{equation}}
\def\endequat{\end{equation}}

    \def\CM{{\mathbb{C}}}

\def\Gb{{\mathbf G}}

\def\Lb{{\mathbf L}}    
    \def\MC{{\mathcal{M}}}

\def\Pb{{\mathbf P}}

\def\O{\Omega}

\def\Sig{\Sigma}

\def\vide{\varnothing}

\def\lexp#1#2{\kern\scriptspace\vphantom{#2}^{#1}\kern-\scriptspace#2}

\def\ge{\hspace{0.1em}\mathop{\geqslant}\nolimits\hspace{0.1em}}

\mathchardef\inferieur="321E
\mathchardef\superieur="321F

\def\eqna{\begin{eqnarray*}}
\def\endeqna{\end{eqnarray*}}

\catcode`\@=11
\long\def\@car#1#2\@nil{#1}
\long\def\@first#1#2{#1}
\long\def\@second#1#2{#2}
\long\def\ifempty#1{\expandafter\ifx\@car#1@\@nil @\@empty
  \expandafter\@first\else\expandafter\@second\fi}
\catcode`\@=12

\def\GL{{\mathrm{GL}}}

\theoremstyle{remark}

\theoremstyle{plain}

\def\BIL{LR}
\def\GAUCHE{L}
\def\CAR{CAR}
\def\FAM{FAM}

\def\xyinj{\ar@{^{(}->}}
\def\xysur{\ar@{->>}}

\bigskip

\makeatletter
\def\hlinewd#1{%
\noalign{\ifnum0=`}\fi\hrule \@height #1 %
\futurelet\reserved@a\@xhline}
\makeatother

\newlength\epaisLigne

\usepackage{multirow,multicol}

\makeatletter
\def\hlinewd#1{%
\noalign{\ifnum0=`}\fi\hrule \@height #1 %
\futurelet\reserved@a\@xhline}
\makeatother

\usepackage{multirow}

\usepackage{lscape}

\def\GL{\operatorname{\Gb\Lb}\nolimits}

\begin{document}

\title{The Maschke octic contains 96 pairwise disjoint lines}

\author{{\sc C\'edric Bonnaf\'e}}
\address{IMAG, Universit\'e de Montpellier, CNRS, Montpellier, France}

\makeatletter
\email{cedric.bonnafe@umontpellier.fr}
\makeatother

\date{\today}

\thanks{The author is partly supported by the ANR 
(Project No ANR-24-CE40-3389, GRAW)}

\begin{abstract}
Miyaoka proved that a smooth surface of degree $d$ in ${\mathbf{P}}^3({\mathbb{C}})$ contains
at most $2d(d-2)$ pairwise disjoint lines. In this note, we verify that the Maschke octic
contains $96$ pairwise disjoint lines, thereby proving that Miyaoka's bound is optimal
for $d=8$.
\end{abstract}

\maketitle
\pagestyle{myheadings}
\markboth{\sc C. Bonnaf\'e}{The Maschke octic contains 96 pairwise disjoint lines}

\bigskip

Let $i$ be a square root of $-1$ in the field $\CM$ of complex numbers.
For better readability, we use dots to represent zeroes in sparse matrices (as is standard
for character tables of finite groups).
Let
$$s_1=
\begin{pmatrix}
. & 1 & . & . \\
1 & . & . & . \\
. & . & 1 & . \\
. & . & . & 1 \\
\end{pmatrix},\qquad s_2=
\begin{pmatrix}
1 & . & . & . \\
. & . & 1 & . \\
. & 1 & . & . \\
. & . & . & 1 \\
\end{pmatrix},\qquad s_3=\begin{pmatrix}
. & -i & . & . \\
i & . & . & . \\
. & . & 1 & . \\
. & . & . & 1
\end{pmatrix},$$
$$s_4=\frac{1}{2}
\begin{pmatrix}
\hphantom{-}1 & -1 & -1 & -1 \\
-1 & \hphantom{-}1 & -1 & -1 \\
-1 & -1 & \hphantom{-}1 & -1 \\
-1 & -1 & -1 & \hphantom{-}1
\end{pmatrix}
\qquad\text{and}\qquad s_5=\begin{pmatrix}
-1 & . & . & . \\
. & 1 & . & . \\
. & . & 1 & . \\
. & . & . & 1
\end{pmatrix}$$
be elements of $\GL_4(\CM)$.
Then $s_1$, $s_2$, $s_3$, $s_4$ and $s_5$ are reflections
of order $2$: they generate a  subgroup of $\Gb\Lb_4(\CM)$ which we denote by $G_{31}$.
It is well-known that $G_{31}$ is finite and is the group denoted by $G_{31}$
in Shephard-Todd classification of finite complex reflection groups~\cite{shephard todd}.
For its natural action on $V=\CM^4$, it is irreducible and primitive.

We denote by $(x,y,z,t)$ the dual basis of the canonical basis of $V$: the algebra $\CM[V]$
of polynomial functions on $V$ is the polynomial algebra $\CM[x,y,z,t]$. If $m$ is a monomial
in $x$, $y$, $z$ and $t$, we denote by $\Sig_4(m)$
the sum of all the monomials obtained from $m$ by permutation of these four variables.
For instance,
$$\Sig_4(xy)=xy+xz+xt+yz+yt+zt\qquad\text{and}\qquad \Sig_4(xyzt)=xyzt.$$
We set
$$f=\Sig_4(x^8)+14\,\Sig_4(x^4y^4)+168\,x^2y^2z^2t^2.$$
Then, $f$ is, up to a scalar, the unique invariant polynomial of $G_{31}$ of degree $8$.
We denote by $\MC$ the surface in $\Pb(V)=\Pb^3(\CM)$ defined as the zero locus of $f$.
It is straightforward to verify with {\sc Magma}~\cite{magma} that $\MC$ is smooth, which implies that the polynomial
$f$ and the surface $\MC$ are irreducible. The surface
$\MC$ is called the {\it Maschke octic}.

It was proved by Boissi\`ere and Sarti that $\MC$ contains $352$ lines~\cite[Theo~3.1]{boissiere}.
It was checked by Sarti and the author~\cite[\S{6.5}]{bonnafe sarti 2} that these $352$ lines are
divided into two $G_{31}$-orbits $\O_{160}$ and $\O_{192}$ of respective cardinality $160$ and $192$.
More precisely, let $L_{160}$ and $L_{192}$ denote the lines in $\Pb^3(\CM)$ defined by
$$
L_{160}:
\begin{cases}
2x + (i + 1)(\sqrt{3} + 1)y=0,\\
2z + (i + 1)(\sqrt{3} + 1)t=0\\
\end{cases}
\qquad \text{and}\qquad
L_{192}:
\begin{cases}
(i+1)(\sqrt{5}+1)x -2(y+z)=0,\\
2(y - z) + (i - 1)(\sqrt{5}+1) t=0.
\end{cases}
$$
Here, $\sqrt{n}$ denotes a square root of $n$ (for $n \in \{3,5\}$).
Then $L_{160}$ and $L_{192}$ are contained
in $\MC$ and $\O_k$ is the orbit of $L_k$ (for $k \in \{160,192\}$).
A computer calculation with {\sc Magma} allows to check easily the main result
of this note:

\bigskip

\noindent{\bf Theorem.} {\it The orbit $\O_{192}$ contains a family of $96$ pairwise
disjoint lines.}

\bigskip

\begin{proof}
Let us give a few details about how this computation is made and some precision on the statement
(again, everything is checked with {\sc Magma}). First,
let
$$
a=\frac{1-i}{2}
\begin{pmatrix}
1 & . & . & 1 \\
. & 1 & 1 & . \\
. & -1 & 1 & . \\
1 & . & . & -1
\end{pmatrix}
\quad\text{and}\quad
b=\begin{pmatrix}
. & . & -1 & . \\
-1 & . & . & . \\
. & 1 & . & . \\
. & . & . & 1
\end{pmatrix}.$$
Then $a$, $b \in G_{31}$. Let $G=\langle a,b\rangle$ and let $\O$ denote the $G$-orbit of the line $L_{192}$.
Then $|G|=1152$ and $|\O|=96$. One can then check that, if $L$, $L' \in \O$ are such that $L \neq L'$,
then $L \cap L'=\vide$. This shows that
one can choose for the family of $96$ pairwise distinct lines in $\O_{192}$
the orbit of $L_{192}$ under the action of $G$.
\end{proof}

\bigskip

\noindent{\it Remarks.} For $d \ge 5$, it is known that there exists a smooth surface of degree $d$
in $\Pb^3(\CM)$ containing $d(d-2)+2$ pairwise disjoint lines~\cite{rams}\footnote{For $d$ odd, this
can be improved to $d(d-2)+4$; see~\cite{rams 5} for $d=5$ and~\cite[Theo.~5.1]{boissiere} for $d \ge 7$.}.
For $d=8$, this gives a smooth octic surface with $50$ pairwise disjoint lines:
the above theorem improves drastically this record.

On the other hand, Miyaoka~\cite[\S{2.2}]{miyaoka} proved that a smooth surface of degree $d \ge 3$ in $\Pb^3(\CM)$
contains at most $2d(d-2)$ pairwise disjoint lines. The above theorem shows that Miyaoka's bound
is optimal for $d=8$.

\bigskip

\noindent{\it Questions.} (1) Miyaoka's bound is optimal for $d\in \{2,3,4,8\}$. Is it
optimal for all $d$?

\medskip

(2) Is the Maschke octic the unique octic surface of $\Pb^3(\CM)$ (up to the action of $\Gb\Lb_4(\CM)$)
containing $96$ pairwise disjoint lines?

\end{document}